\documentclass[12pt]{amsart} 
\usepackage{amssymb}
\usepackage{latexsym}
\usepackage{amsthm}
\usepackage{graphicx}
\usepackage{epstopdf}
\usepackage[all]{xy}
\usepackage{verbatim}
\input epsf

\newtheorem{prop}{Proposition}

\newtheorem{definition}{Definition}
\newtheorem{thm}{Theorem}

\newtheorem{cor}{Corollary}
\newtheorem{lemma}{Lemma}

\theoremstyle{definition}
\newcommand{\lk}{\mathrm{lk}}

\title{Intrinsic linking and knotting are arbitrarily complex}
\author{Erica Flapan, Blake Mellor, and Ramin Naimi}



\address{Department of Mathematics, Pomona College, Claremont, CA 91711, USA}
\email{eflapan@pomona.edu}

\address{Department of Mathematics, Loyola Marymount University, Los 
Angeles, CA 90045, USA}
\email{bmellor@lmu.edu}

\address{Department of Mathematics, Occidental College, Los Angeles, 
CA 90041, USA}
\email{rnaimi@oxy.edu}

\begin{document}

\date \today
\maketitle
\begin{abstract}
We show that,
given any $n$ and $\alpha$,
every embedding of any sufficiently large complete graph in $\mathbb{R}^3$
contains
an oriented link with components $Q_1$, \dots, $Q_n$
such that for every $i\not =j$, $|\lk(Q_i,Q_j)|\geq\alpha$
and $|a_2(Q_i)|\geq\alpha$, where $a_{2}(Q_i)$
denotes the second  coefficient
of the Conway polynomial of $Q_i$.
  \end{abstract}

\section{Introduction}
The study of embeddings of graphs in $\mathbb{R}^3$ is a natural 
extension of knot theory.  However, in contrast with knots whose 
properties depend only on their extrinsic topology, there is a rich 
interplay between the intrinsic structure of a graph and the 
extrinsic topology of all embeddings of the graph in $\mathbb{R}^3$. 
Conway and Gordon \cite{cg} obtained groundbreaking results of this 
nature by showing that every embedding of the complete graph $K_6$ in 
$\mathbb{R}^3$ contains a non-trivial link and every embedding of 
$K_7$ in $\mathbb{R}^3$ contains a non-trivial knot. Because this 
type of linking and knotting is intrinsic to the graph itself rather 
than depending on the particular embedding of the graph in 
$\mathbb{R}^3$,  $K_6$ is said to be {\it intrinsically linked} and 
$K_7$ is said to be {\it intrinsically knotted}.  On the other hand, 
Conway and Gordon \cite{cg} illustrated an embedding of $K_{6}$ such
that the only non-trivial link $L_{1}\cup L_{2}$
contained in $K_{6}$ is the Hopf link (which has
$ \vert lk(L_{1},L_{2}) \vert =1$); and they
illustrated an embedding of $K_{7}$ such that the
only non-trivial knot $Q$ contained in $K_{7}$ is
the trefoil knot (which has $ \vert a_{2}(Q)
\vert =1$, where $a_{2}(Q)$
denotes the second  coefficient
of the Conway polynomial of $Q$).  In this sense, we see that $K_{6}$ 
exhibits the
simplest type of intrinsic linking and $K_{7}$
exhibits the simplest type of intrinsic knotting.

More recently, it has been shown that for sufficiently large
values of $r$, the complete graph $K_{r}$
exhibits more complex types of intrinsic linking
and knotting.  In particular, Flapan  \cite{fl} showed that for
every $\lambda\in {\mathbb N}$, there
is a complete graph $K_r$ such
that every
embedding of
$K_{r}$ in ${\mathbb  R}^{3}$
contains both a 2-component oriented link $L$
whose  linking number is at least
$\lambda$ and a knot
$Q$ with $ \vert a_{2}(Q) \vert
\geq  \lambda$ (though $L$ and $Q$ have no particular relationship). 
Fleming \cite{ft} shows that for any $n\in \mathbb{N}$, there is a 
graph $G$ such that every embedding of $G$ in $\mathbb{R}^3$ contains 
a non-split link of $n+1$ components where $n$ of the components are 
non-trivial knots.  In the current paper, we show that for 
sufficiently large complete graphs intrinsic linking with knotted 
components is arbitrarily complex both in terms of linking number and 
in terms of the knotting of every component.  In particular, our main 
result is the following.

\setcounter{thm}{1}
\begin{thm}\label{main thm}  For every $n$, $\alpha\in \mathbb{N}$, 
there is a complete graph $K_r$ such that every embedding of $K_r$ in 
$\mathbb{R}^3$ contains
an oriented link with components $Q_1$, \dots, $Q_n$
such that for every $i \not =j$, $|\lk(Q_i,Q_j)|\geq\alpha$ and 
$|a_2(Q_i)|\geq\alpha$.
\end{thm}

If linking is measured with linking number and knotting is measured with $a_2$, then this is the strongest result one could hope for about the complexity of simultaneous intrinsic knotting and linking.  Furthermore, observe that for a given $c\in \mathbb{N}$, there are only finitely many knots whose minimal crossing number is less than or equal to $c$.  If we pick $\lambda$ larger than the $|a_2|$ of all of the knots with minimal crossing number less than or equal to $c$, then the knots $Q_1$, \dots, $Q_n$ given by Theorem \ref{main thm} will each have minimal crossing number greater than $c$.  It follows that the complexity of intrinsic knotting as measured by the crossing number can also be made arbitrarily large.

In order to prove our main result, we first prove in Section 2 that 
intrinsic linking is arbitrarily complex in the sense of the 
structure of a link.  In particular, we prove the following.

\setcounter{thm}{0}
  \begin{thm}\label{T:complete}  For every $n$, $\lambda\in 
\mathbb{N}$, there is a complete graph $K_r$ such that every 
embedding of $K_r$ in $\mathbb{R}^3$ contains
an oriented link with components $Q_1$, \dots, $Q_n$
such that for every $i\not =j$, $|\lk(Q_i,Q_j)|\geq\lambda$.
\end{thm}

The idea of the proof of Theorem \ref{main thm} is then as follows. 
Taniyama and Yasuhara \cite{ty} have shown that for any embedding of 
the pseudo-graph $D_4$ (see Figure \ref{collapseD4}), the product of 
the linking numbers of the two opposite pairs of cycles, $C_i$ and 
$C_j$, is related to the sum of the $a_2$'s of all of the Hamiltonian 
cycles $Q$ according to the formula:

\begin{equation*}  \sum _{Q \in
S  }|a_{2}(Q )|\geq |\lk(C_1,C_3)\lk(C_2,C_4)|
\end{equation*}

We use Theorem \ref{T:complete} to find a complete graph $K_r$, so 
that every embedding of $K_r$ contains a link with a  ``large 
enough'' number of components all pairwise linked with linking number 
at least $\lambda$. We then use a recursive argument to successively 
exchange opposite pairs of linked cycles in a $D_4$ (which is a minor 
of $K_r$) for a knotted Hamiltonian cycle in $D_4$ which is linked 
with all of the previous knotted cycles in the construction.  We do 
this in such a way that for all of the knotted cycles $Q_i$ and $Q_j$ 
we have $|\lk(Q_i,Q_j)|\geq\lambda\geq \alpha$ and 
$|a_2(Q_i)|\geq\lambda^2/16\geq \alpha$.
\bigskip

We wish to thank Alan Tarr for suggesting that we represent the 
structure of a link with a linking pattern, for reading a preliminary 
draft of this paper, and for suggesting that we prove the current 
version of Theorem \ref{T:complete}.  We also want to thank Kouki Taniyama for asking us whether the conclusion of Theorem 1 could be strengthened by removing the absolute value on the linking number.  The current version of Proposition \ref{P:nostrong}, together with Corollary \ref{C:3link}, and Proposition \ref{P:Ramsey} grew out of our response to Taniyama's question.

\section{Intrinsic linking}

We will use a weighted graph to describe the structure of a link as follows.

\medskip

\begin{definition}
Given an oriented link $L$ with components $L_1$, \dots, $L_n$, the 
\emph{linking pattern} of $L$ is the graph with vertices $v_1$, 
\dots, $v_n$, such that there is an edge between $v_i$ and $v_j$ if 
and only if $\mathrm{lk}(L_i,L_j) \neq 0$.  The \emph{weighted 
linking pattern} of $L$ is the linking pattern with a  \emph{weight} 
assigned to each edge $\{v_i,v_j\}$ representing the value of 
$|\mathrm{lk}(L_i,L_j)|$.
\end{definition}
\medskip

For example, the linking pattern of a Hopf link is a single edge.  As 
another example, consider a {\it keyring} link.  That is, a link 
consisting of a {\it ring}, $J$, and {\it keys},  $L_1$, \dots, 
$L_n$,  such that $\mathrm{lk}(J, L_i)\not=0$ and 
$\mathrm{lk}(L_i,L_j)=0$ for all $i\not =j$.  The linking pattern of 
a keyring link is an {$n$-star} (i.e., a graph consisting of $n$ 
vertices all connected to a single additional vertex).

Many results about intrinsic linking use the mod 2 linking number, 
$\omega(J,L)\linebreak
=\mathrm{lk}(J,L) \mod 2$ as a simpler measure of linking 
than the ordinary linking number.  Thus we will also use the 
following definition.

\medskip
\begin{definition}
Given a link $L$ with components $L_1$, \dots, $L_n$, the \emph{mod 2 
linking pattern} of $L$ is the graph with vertices $v_1$, \dots, 
$v_n$, such that there is an edge between $v_i$ and $v_j$  if and 
only if $\omega(L_i,L_j)=1$.
\end{definition}
\medskip

Using this terminology, Conway and Gordon \cite{cg} and Sachs \cite 
{S1}\cite{S2} showed that every embedding of $K_6$ in $\mathbb{R}^3$ 
contains a link whose mod 2 linking pattern is a single edge. 
Fleming and Diesl \cite{fd} showed that there is a graph $G$ such 
that every embedding of $G$ in $\mathbb{R}^3$ contains a link $J \cup 
L_1 \cup \dots \cup L_n$  where $\omega(J, L_i)=1$ for each $i$. 
Thus every embedding of $G$ in $\mathbb{R}^3$ contains a link whose 
mod 2 linking pattern contains an $n$-star, possibly with additional 
edges.  We call such a link a {\it generalized keyring} link, since 
some of the $L_i$'s may be linked with one another.  

Prior to this paper, a chain of $n$ edges and a circle of $n$ edges  \cite{ffnp}, and an $n$-star  \cite{fd}
were the only linking patterns $\Gamma$ which were known to have the 
property that for some graph $G$ every embedding of $G$ in 
$\mathbb{R}^3$ contains a link whose linking pattern contains 
$\Gamma$.  We prove in Theorem \ref{T:complete} that for every 
complete graph $K_n$, there is a graph $G$ such that every embedding 
of $G$ in $\mathbb{R}^3$ contains a link whose linking pattern is 
$K_n$.  It follows that for any linking pattern $\Gamma$, there is a 
graph $G$ such that every embedding of $G$ in $\mathbb{R}^3$ contains 
a link whose linking pattern contains $\Gamma$. In particular, using 
the language of linking patterns we prove the following restatement 
of Theorem \ref{T:complete}.
\medskip

\setcounter{thm}{0}
\begin{thm} \label{T:complete}
  Let $\lambda \in \mathbb{N}$.  For every $n\in \mathbb{N}$, there is 
a graph $G$ such that every embedding of $G$ in $\mathbb{R}^3$ 
contains a link whose linking pattern is $K_n$ with every weight at 
least $\lambda$.\end{thm}

This theorem implies that given any $n$, $\lambda\in \mathbb{N}$, 
there exists a graph $G$ which has the property that every embedding 
of $G$ contains an $n$-component link all of whose components are 
pairwise linked with the absolute value of their linking number at least $\lambda$ (this was our 
statement of Theorem \ref{T:complete} in the introduction).
Furthermore, in Proposition \ref{P:nostrong} we show that complete 
graphs are the only linking patterns that have the property described 
by Theorem \ref{T:complete}.  In other words, complete graphs can be 
said to be the only {\it intrinsic} linking patterns.

  Before we prove Theorem \ref{T:complete}, we will show in 
Proposition \ref{T:bipartite} that for every $n$, there is a graph 
such that every embedding of the graph contains a link whose linking 
pattern contains $K_{n,n}$.  We will use this result to prove Theorem 
\ref{T:complete}.  In fact,  in Proposition \ref{T:bipartite}, we 
prove the stronger result that we can ensure that all of the weights 
are odd.  It is  an open question whether this stronger formulation 
can be extended to the linking pattern $K_n$.

  In order to prove Proposition \ref{T:bipartite}, we need the 
following lemma which allows us to combine many pairs of linked 
cycles into a single cycle that links some proportion of the original 
components.  Throughout the paper we use the term {\it cycle} to mean 
a simple closed curve within a graph and the notation $J\triangledown 
L=\overline{(J\cup L)-(J\cap L)}$ for the closure of the symmetric 
difference.

\medskip

\begin{lemma} \label{L:multimerge}
Let $K_p$ be embedded in $\mathbb{R}^3$ such that it contains a  link 
with components $J_1$, \dots, $J_{n^2}$ and $X_1$, \dots, $X_{n^2}$, 
and $\omega(J_i,X_i)=1$ for every $i\leq n^2$.  Then there is a cycle 
$Z$ in $K_p$ with vertices on $J_1\cup \dots\cup J_{n^2}$, and an 
index set $I$ with $|I|\geq \frac{n}{2}$, such that $\omega(Z,X_j)=1$ 
for all $j\in I$.
\end{lemma}

\begin{figure}[h]
\includegraphics{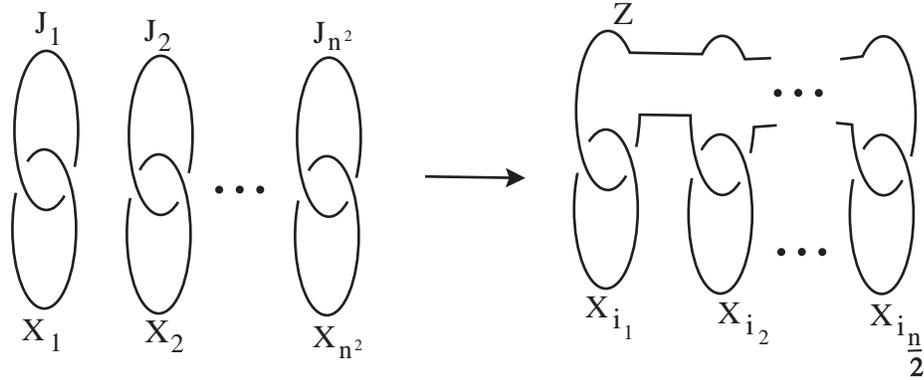}
\caption{Illustration of Lemma 1}
\label{lemma1}
\end{figure}

\begin{proof}  We begin by creating a cycle $C$ in $K_p$ which 
cyclicly joins the $J_i$ as follows.  For each $i=1$, \dots, 
${n^2}$, let $u_i$ and $w_i$ be vertices on the cycle $J_i$, and let 
$q_i$ be a path on $J_i$ from $u_i$ to $w_i$.  For $i=1$, \dots, 
$n^2-1$, let $e_i$ be the edge in $K_p$ from $w_i$ to $u_{i+1}$, and 
let $e_{n^2}$ be the edge in $K_p$ from $w_{n^2}$ to $u_1$.  Let $C$ 
be the cycle $\bigcup_{i=1}^{n^2}{q_i\cup e_i}$.  Observe that if at 
least $\frac{n}{2}$ of the $\omega(C,X_j)$'s are equal to $1$, then 
we are done by letting $Z = C$.  So we shall assume that fewer than 
$\frac{n}{2}$ of the $\omega(C,X_j)$'s are equal to $1$.

We create an $n^2\times n^2$ matrix $M$ with entries in 
$\mathbb{Z}_2$ as follows.  For each $i,j=1$, \dots, $n^2$, let the 
entry $M_{ij}=\omega(J_i,X_j)$.  By hypothesis, for each $i$, 
$M_{ii}=1$.  Using row reduction mod 2, let $M'$ denote the reduced 
row-echelon form of $M$.  Observe that since every column of $M$ 
contains a $1$, every column of $M'$ also contains a $1$.

Next we will add together rows of $M'$ as follows to create a row 
vector $V$, at least $n$ of whose entries are $1$'s. Let $r$ denote 
the rank of $M'$ over $\mathbb{Z}_2$.  First suppose that $r\geq n$. 
In this case, let $V$ be the vector obtained by adding together all 
of the non-trivial rows of $M'$ modulo 2.  Then $V$ has at least 
$r\geq n$ entries which are $1$'s. So $V$ is the desired vector.  On 
the other hand, suppose that $r<n$.  Then $\frac{n^2}{r}> n$. 
Observe that $M'$ contains $r$ non-trivial rows and has $n^2$ columns 
which contain 1's.  So by the Pigeonhole Principle, some row of $M'$ 
has at least $\frac{n^2}{r} > n$  entries which are 1's. In this 
case, let $V$ be the vector representing this row.

In either case, $V$ can be written as the sum of some of the rows of 
$M$, say rows $i_1$, \dots, $i_k$.   Thus, for each $j=1$, \dots, 
$n^2$, the $j^{th}$ entry of the vector $V$ is 
$V_j=\omega(J_{i_1},X_j)+\dots+\omega(J_{i_k},X_j) \mod{2}$.  Recall 
that $V$ was chosen so that at least $n$ of the $V_j$'s are equal to 
$1$.  Also we assumed that fewer than $\frac{n}{2}$ of the 
$\omega(C,X_j)$'s are equal to $1$.  Let $I$ denote the subset of 
$\{1,\dots, n^2\}$ such that for each $j\in I$, we have 
simultaneously $V_j=1$ and $\omega(C,X_j)=0$.  Then $|I|> 
n-\frac{n}{2}=\frac{n}{2}$.  Now let $Z = C \triangledown J_{i_1} 
\triangledown \dots \triangledown J_{i_k}$.  Then for each $j\in I$, 
$\omega(Z,X_j)=1$, as required.\end{proof}

\medskip

Recall the following definition from graph theory.  The {\it complete 
bipartite graph}, $K_{m,n}$, is defined as the graph whose vertices 
are partitioned into subsets $P_1$ and $P_2$, where $P_1$ contains 
$m$ vertices, $P_2$ contains $n$ vertices, and there is an edge 
between two vertices if and only if one vertex is in $P_1$ and the 
other is in $P_2$.
\medskip

\begin{prop} \label{T:bipartite}
For every $n$, there is a graph $G$ such that every embedding of $G$ 
in $\mathbb{R}^3$ contains a link whose mod 2 linking pattern 
contains the complete bipartite graph $K_{n,n}$.
\end{prop}

\begin{proof}   Let $n$ be given.  It was shown in Fleming and Diesl 
\cite{fd} Lemma 2.3, that there exists some $p$ such 
that every embedding of $K_p$ in $\mathbb{R}^3$ contains a mod 2 
generalized keyring link with a ring and $n$ keys.  Observe that the existence of such a $p$ also follows from our Lemma \ref{L:multimerge} together with Conway and Gordon's \cite{cg} result that every embedding of $K_6$ contains a link $J\cup X$ such that $\omega(J,X)=1$.  

Let 
$m=\frac{(4n)^{2^n}}{4}$.  Every embedding of $K_{mp}$ in 
$\mathbb{R}^3$ contains $m$ disjoint mod 2 generalized keyring links 
each with a ring and $n$ keys.  We will prove that every embedding of 
$K_{mp}$ in $\mathbb{R}^3$ contains a link whose mod 2 linking 
pattern contains $K_{n,n}$.

Let $K_{mp}$ be embedded in $\mathbb{R}^3$.  Let $X_1$, \dots, $X_m$ 
denote the rings of the generalized keyring links in the $m$ disjoint 
copies of $K_p$ in $K_{mp}$.  For each $i\leq m$ and $j\leq n$, let 
$J_{ij}$ be a key on the ring $X_i$.

Since $\omega(J_{i1},X_i)=1$ for all $i=1\leq m$, we can apply Lemma 
\ref{L:multimerge} to the link in $K_{mp}$ with components $J_{11}$, 
$J_{21}$, \dots, $J_{m1}$ and $X_1$, \dots, $X_m$.  This gives us a 
cycle $Z_1$ with vertices on $J_{11}\cup \dots\cup J_{m1}$ and an 
index set $I_1$ with $|I_1|\geq 
\frac{\sqrt{m}}{2}=\frac{(4n)^{2^{n-1}}}{4}=r_1$, such that for each 
$i\in I_1$, $\omega(Z_1,X_i)=1$.  Now since $\omega(J_{i2},X_i)=1$ 
for all $i\in I_1$, we can apply Lemma \ref{L:multimerge} to the link 
in $K_{mp}$ whose components are all those $J_{i2}$ and $X_i$ with 
$i\in I_1$.  This gives us a cycle $Z_2$ with vertices on $J_{12}\cup 
\dots\cup J_{m2}$ and an index set $I_2\subseteq I_1$ with $|I_2|\geq 
\frac{\sqrt{r_1}}{2}=\frac{(4n)^{2^{n-2}}}{4}=r_2$, such that for 
each $i\in I_2$, $\omega(Z_2,X_i)=1$. Continue this process to get 
disjoint cycles $Z_1$, \dots, $Z_n$ in $K_{mp}$ and index sets 
$I_n\subseteq \dots \subseteq I_1$ with $|I_n|\geq 
\frac{\sqrt{r_{n-1}}}{2}=\frac{(4n)^{2^{n-n}}}{4}=n$ such that for 
every $i\in I_n$ and every $j\leq n$, $\omega(Z_j,X_i)=1$.

Thus the mod 2 linking pattern of the link with components $Z_1$, 
\dots $Z_n$ and all those $X_i$ with $i\in I_n$ contains $K_{n,n}$. 
Hence  every embedding of $K_{mp}$ in $\mathbb{R}^3$ contains a link 
whose mod 2 linking pattern contains $K_{n,n}$.
\end{proof}
\medskip

Before we prove Theorem \ref{T:complete}, we need one more lemma that 
allows us to combine components of distinct links.

\medskip

\begin{lemma} \label{L:multijoin}
Let $\lambda \in \mathbb{N}$.  Let $K_p$ be embedded in 
$\mathbb{R}^3$ such that it contains a link with oriented components 
$J_1$, \dots, $J_r$, $L_1$, \dots, $L_q$, $X_1$, \dots, $X_m$, and 
$Y_1$, \dots, $Y_n$ where $r \geq m(2\lambda +1)2^m$ and $q \geq 
(m+n)(2\lambda +1)3^m2^n$, and for every $i, j,\alpha, \beta$, 
$\mathrm{lk}(J_i,X_\alpha) \neq 0$ and $\mathrm{lk}(L_j, Y_\beta) 
\neq 0$.  Then $K_p$ contains a cycle $Z$ with vertices in 
$J_1\cup\dots\cup J_r\cup L_1\cup\dots\cup L_q$ such that for every 
$\alpha$ and $\beta$, $|\mathrm{lk}(Z, X_\alpha)| >\lambda$ and 
$|\mathrm{lk}(Z,Y_\beta)|>\lambda$.
\end{lemma}

\begin{figure}[h]
\includegraphics[width=5in]{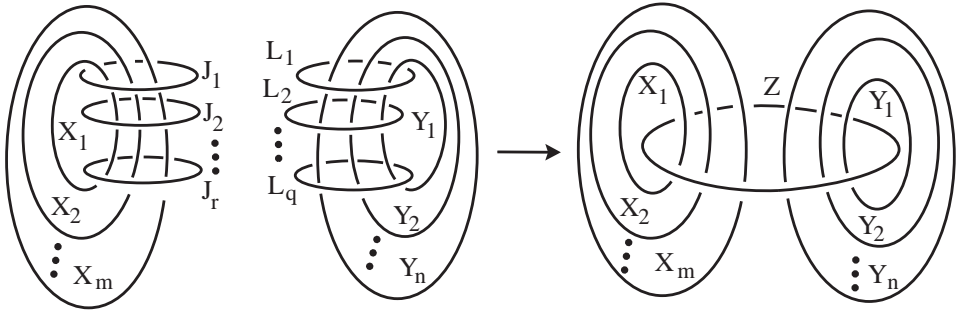}
\caption{Illustration of Lemma 2}
\label{lemma2}
\end{figure}

\begin{proof}  First we consider the signs of the linking numbers 
$\mathrm{lk}(J_i,X_1)$, for $i\leq r$.   At least half of these 
linking numbers have the same sign.  Thus without loss of generality, 
$J_1$, \dots, $ J_{\lceil{\frac{r}{2}}\rceil}$  all have the same 
sign linking number with $X_1$.  Furthermore, we may assume these 
linking numbers are all positive (otherwise, simply reverse the 
orientation of $X_1$).  Now, consider the signs of the linking 
numbers $\mathrm{lk}(J_i,X_2)$, for $i\leq\lceil{\frac{r}{2}}\rceil$. 
At least half of these have the same sign.  We continue this process 
for each subsequent $X_\alpha$.  In this way, we end up with cycles 
$J_1$, \dots, $J_{r'}$ which each have positive linking number with 
every $X_\alpha$, and $r' \geq \frac{r}{2^m} \geq m(2\lambda +1)$. 
 From now on the only $J_i$ that we consider will be $J_1$, \dots, $ 
J_{m(2\lambda +1)}$, which all have positive linking number with 
every $X_\alpha$.

Next we consider the signs of the linking numbers 
$\mathrm{lk}(L_j,Y_1)$, for $j\leq  q$.  By using an argument 
analogous to the one above, we end up with cycles $L_1$, \dots, 
$L_{q'}$ which each have positive linking number with every 
$Y_\beta$, and $q' \geq \frac{q}{2^n} \geq (m+n)(2\lambda +1)3^m$.

Now we consider the signs (positive, negative, or zero) of the 
linking numbers $\mathrm{lk}(L_j,X_1)$, for $j\leq q'$.  By the 
Pigeonhole Principle at least one third of these signs are the same. 
Without loss of generality, we can assume that $L_1$, \dots, 
$L_{\lceil{\frac{q'}{3}}\rceil}$ each have the same sign (positive, 
negative, or zero) linking number with $X_1$.  We continue this 
process for each subsequent $X_\alpha$.  In this way, we end up with 
cycles $L_1$, \dots, $L_{q''}$ which each have the same sign 
(positive, negative, or zero) linking number with every $X_\alpha$, 
and $q'' \geq \frac{q'}{3^m} \geq (m+n)(2\lambda +1)$.  From now on 
the only $L_j$ that we consider are $L_1$, \dots, $ L_{(m+n)(2\lambda 
+1)}$, which all have positive linking number with every $Y_\beta$ 
and all have the same sign linking number with every $X_\alpha$.

We shall create a cycle $C_0$ which cyclicly joins $J_1$, \dots, 
$J_{m(2\lambda +1)}$, and $L_1$, \dots, $ L_{(m+n)(2\lambda +1)}$  as 
follows.  For $i\leq m(2\lambda +1)$, let $u_i$ and $w_i$ be vertices 
on $J_i$, and for $j\leq (m+n)(2\lambda +1)$, let $u_{m(2\lambda 
+1)+j}$ and $w_{m(2\lambda +1)+j}$ be vertices on $L_j$.  Also for 
$k\leq (2m+n)(2\lambda +1)-1$, let $e_k$ be the edge in $K_p$ from 
$w_k$ to $u_{k+1}$, and let $e_{(2m+n)(2\lambda +1)}$ be the edge in 
$K_p$ from $w_{(2m+n)(2\lambda +1)}$ to $u_1$.  Finally, for $i\leq 
m(2\lambda +1)$, let $q_i$ be the path on $J_i$ from $u_i$ to $w_i$ 
which travels {\it opposite} to the orientation on $J_i$,  and for 
$j\leq (m+n)(2\lambda +1)$, let  $q_{m(2\lambda +1)+j}$  be the path 
on $L_{j}$ from $u_{m(2\lambda +1)+j}$ to $w_{m(2\lambda +1)+j}$ 
which travels {\it opposite} to the orientation of $L_j$.  Now let 
$C_0 = \bigcup_{k=1}^{(2m+n)(2\lambda +1)}{e_k \cup q_k}$ with the 
induced orientation.
Also for each $s\leq m(2\lambda +1)$, let $C_s = C_0 \triangledown 
J_1 \triangledown \dots \triangledown J_s$.

  Observe that for every $\alpha$ and $s\leq m(2\lambda+1)-1$, 
$\mathrm{lk}(C_{s+1},X_\alpha) 
=\mathrm{lk}(C_s,X_\alpha)+\mathrm{lk}(J_{s+1},X_\alpha)$ and 
$\mathrm{lk}(J_{s+1},X_\alpha) > 0$.  Hence for a given $\alpha$, 
$\mathrm{lk}(C_s,X_{\alpha})$ is a strictly increasing function of 
$s$.  In particular, for a given $\alpha$, the values of 
$\mathrm{lk}(C_s,X_{\alpha})$ are distinct for different values of 
$s$.  As there are $2\lambda +1$ distinct values of $a$ with $|a|\leq 
\lambda$, for a given $\alpha$ there are at most $2\lambda +1$ values 
of $s$ such that $|\mathrm{lk}(C_s,X_\alpha)|\leq \lambda$.  Since 
there are $m$ values of $\alpha$, there are at most $m(2\lambda +1)$ 
values of $s$ such that there is an $\alpha$ with 
$|\mathrm{lk}(C_s,X_\alpha)|\leq \lambda$.  Now by the Pigeonhole 
Principle, since there are $m(2\lambda +1)+1$ values of $s$ 
(including $s=0$),  there must be at least one $C_s$ such that 
$|\mathrm{lk}(C_s,X_\alpha)|>\lambda$ for every $\alpha$.  Let $D_0$ 
denote such a $C_s$.

Recall that each $Y_\beta$ has positive linking number with $L_1$, 
\dots, \linebreak $L_{(m+n)(2\lambda +1)}$, and each $X_\alpha$ has 
the same sign (positive, negative or zero) linking number with $L_1$, 
\dots, $L_{(m+n)(2\lambda +1)}$.   In fact, by changing the 
orientation of some $X_\alpha$, we can assume that each $X_\alpha$ 
has nonnegative linking number with $L_1$, \dots, $L_{(m+n)(2\lambda 
+1)}$.  Note that changing the orientation of a particular $X_\alpha$ 
does not change the fact that $|\mathrm{lk}(D_0,X_\alpha)|>\lambda$. 
Now let $S$ denote the set of all the $Y_\beta$'s together with those 
$X_\alpha$'s which have positive linking number with all of $L_1$, 
\dots, $L_{(m+n)(2\lambda +1)}$.  For each $t\leq (m+n)(2\lambda 
+1)$, let $D_t = D_0 \triangledown L_1 \triangledown \dots 
\triangledown L_t$.

Since $\alpha \leq m$ and $\beta \leq n$,  the set $S$ contains at 
most $m+n$ cycles.  Also, each cycle in $S$ has positive linking 
number with $L_1$, \dots, $L_{(m+n)(2\lambda +1)}$.  Thus we can use 
the same argument as the one we used for the $C_s$ to show that there 
is some $D_t$ such that $|\mathrm{lk}(D_t,A)|>\lambda$ for each $A\in 
S$.  Let $Z$ denote such a $D_t$.  Finally, observe that for each 
$X_\alpha$ not in $S$, $\mathrm{lk}(L_j,X_\alpha) = 0$ for all $j\leq 
(m+n)(2\lambda +1)$.  Hence if $X_\alpha\not \in S$, then 
$|\mathrm{lk}(Z, X_\alpha)| = |\mathrm{lk}(D_0, X_\alpha)| >\lambda$. 
So, for every $\alpha$ and $\beta$, $|\mathrm{lk}(Z, 
X_\alpha)| >\lambda$ and $|\mathrm{lk}(Z,Y_\beta)|>\lambda$ as 
desired.
\end{proof}

\medskip

In the proof of Theorem \ref{T:complete}, we will use Lemma 
\ref{L:multijoin} together with the following definition from graph 
theory.  A {\it complete $m$-partite graph}, is defined as a graph 
whose vertices are partitioned into subsets $P_1$, \dots, $P_m$, and 
there is an edge between two vertices if and only if the two vertices 
are in distinct subsets of the partition.

\medskip
\setcounter{thm}{0}
\begin{thm} \label{T:complete}
  Let $\lambda \in \mathbb{N}$.  For every $n\in \mathbb{N}$, there is 
a graph $G$ such that every embedding of $G$ in $\mathbb{R}^3$ 
contains a link whose linking pattern is $K_n$ with every weight at 
least $\lambda$.\end{thm}

\begin{proof}  For each $m$, $n\in \mathbb{N}$, let $H(n,m)$ denote 
the complete $(n+2)$-partite graph with partitions $P_1$ and $P_2$ 
containing $m$ vertices each and partitions $Q_1$, \dots,$Q_n$ each 
containing a single vertex. We will prove by induction on $n$ that 
for each $n\geq 0$, for every $m\geq 1$, there is a graph $G$ such 
that every embedding of $G$ in $\mathbb{R}^3$ contains a link whose 
weighted linking pattern contains $H(n,m)$, and the weight of every 
edge between vertices in $Q_i$ and $Q_j$ is greater than $\lambda$.

The base case is for $n=0$.  Observe that $H(0,m)$ is the complete 
bipartite graph $K_{m,m}$.  Hence by Proposition \ref{T:bipartite} we 
know that for every $m$, there is a graph $G$ such that every 
embedding of $G$ in $\mathbb{R}^3$ contains a link whose linking 
pattern contains $H(0,m)$.

As our induction hypothesis we suppose that for some $n\geq 0$, for 
every $m\geq 1$ there is a graph $G$ such that every embedding of $G$ 
in $\mathbb{R}^3$ contains a link whose linking pattern contains 
$H(n,m)$, and the weight of every edge between vertices in $Q_i$ and 
$Q_j$ is greater than $\lambda$.  Let $m$ be given.  Let 
$q=(2m+n)(2\lambda +1)3^m2^{m+n}$ and let $s=m+q$.  Consider the 
graph $H(n,s)$, with partitions $P_1$, $P_2$, $Q_1$,\dots, $Q_{n}$, 
where the vertices in each partition are denoted as follows.  The $s$ 
vertices in $P_1$ are denoted by $X_1$, \dots, $X_m$, $L_1$, \dots, 
$L_q$.  The $s$ vertices in $P_2$ are denoted by $Y_{1}$, \dots, 
$Y_{m}$, $J_1$, \dots, $J_q$.  For each $i\leq n$, the partition 
$Q_i$ contains a single vertex denoted by $Y_{m+i}$.

It follows from our inductive hypothesis that there is a graph $G$ 
such that every embedding of $G$ in $\mathbb{R}^3$ contains a link 
$L$ whose linking pattern contains $H(n,s)$, and for every $i\not = 
j$ the weight of the edge between vertices $Y_{m+i}$ and $Y_{m+j}$ is 
greater than $\lambda$. Without loss of generality, $G$ is a complete 
graph $K_p$.  Let $K_p$ be embedded in $\mathbb{R}^3$. We will prove 
that $K_p$ also contains a link whose weighted linking pattern 
contains $H(n+1, m)$ with the desired weights.   We shall abuse 
notation and let each of the components of the link $L$ in $K_p$ be 
denoted by the name of the vertex that represents that component in 
the linking pattern described above.

We can apply Lemma \ref{L:multijoin} to the oriented link with 
components $J_{ 1}$, \dots, $ J_{ q}$, $L_ { 1}$, \dots, $ L_ q$, 
$X_1$, \dots, $ X_m$, and $Y_{1}$, \dots, $Y_{m+n}$ in $K_p$ with 
$r=q=(2m+n)(2\lambda +1)3^m2^{m+n}$, to get a cycle $Y_{m+n+1}$ with 
vertices in $J_{ 1}\cup \dots\cup J_{ q}\cup L_ { 1}\cup \dots\cup 
L_{ q}$, such that for every $\alpha\leq n$ and $\beta\leq m+n$, 
$|\mathrm{lk}(Y_{m+n+1}, X_\alpha)| >\lambda$ and 
$|\mathrm{lk}(Y_{m+n+1},Y_\beta)|>\lambda$.

Thus $K_p$ contains a link $L'$ with components $X_1$, \dots, $X_m$, 
$Y_1$, \dots, $Y_{m+n+1}$.  The components of $L'$ can be partitioned 
into subsets $P'_1$, $P'_2$, $Q'_1$, \dots, $ Q'_{n+1}$, where $P'_1$ 
contains  $X_1$,\dots, $X_m$, $P'_2$ contains $Y_1$, \dots, $Y_m$, 
and for each $i=1$, \dots, $n+1$, $Q'_i$ contains $Y_{m+i}$. 
Furthermore every component in one partition is linked with every 
component in all the other partitions.  Also, for each $i$, $j\leq 
n+1$ with $i\not = j$, $|\mathrm{lk}(Y_{m+i},Y_{m+j})|>\lambda$.  It 
follows that the weighted linking pattern of $L'$ contains 
$H(n+1,m)$, and the weight of every edge between vertices $Q'_i$ and 
$Q'_j$ is greater than $\lambda$.

Thus, we have shown that for every $n\geq 0$ and $m\geq 1$, there is 
a graph $G$ such that every embedding of $G$ in $\mathbb{R}^3$ 
contains a link whose weighted linking pattern contains $H(n,m)$, and 
the weight of every edge between vertices $Y_{m+i}$ and $Y_{m+j}$ in 
$Q_i$ and $Q_j$ respectively is greater than $\lambda$.  Observe that 
the subgraph of $H(n,m)$ consisting of the vertices $Y_{m+1}$, \dots, 
$Y_{m+n}$ in $Q_1$, \dots, $Q_n$ respectively, together with the 
edges between them is the complete graph $K_n$.  Hence every 
embedding of $G$ in $\mathbb{R}^3$ contains a link whose linking 
pattern is $K_n$, and the weight of every edge of $K_n$ is greater 
than $\lambda$.
\end{proof}

\medskip

    The following proposition shows that every graph
$G$ has some embedding in $\mathbb{R}^3$ such that the linking
pattern of every link in that embedding is a complete graph.  Hence complete graphs are the only
linking patterns which have the property described by Theorem
\ref{T:complete}.   

\medskip

\begin{prop}\label{P:nostrong}
Given any $\lambda > 0$ and
any graph $G$, there exists an embedding $G'$ of $G$ in $\mathbb{R}^3$ and an orientation of the cycles in $G'$ such that
for every pair of disjoint oriented cycles $C_1$ and $C_2$ in $G'$,
 $\lk(C_1, C_2) \ge \lambda$.
\end{prop}

\begin{proof}
We start with any embedding of the graph $G$ in $\mathbb{R}^3$.
We arbitrarily assign to each edge a unique positive integer which will be the weight of that edge.  The weights give the set of edges of $G$ a linear ordering.  We put an arbitrary orientation on each edge,
and then orient each cycle according to the orientation of its edge whose weight is the largest. For each pair of disjoint edges $(e_i, f_i)$  we denote the smaller weighted edge by $e_i$ and the larger weighted edge by $f_i$.  Note that a given edge may occur in multiple pairs, sometimes as the smaller weighted edge and sometimes as the larger weighted edge.
Let $(e_0, f_0), \cdots, (e_n,f_n)$
denote the set of all pairs of disjoint edges
ordered lexicographically from smallest weighted pair to largest weighted pair.
We describe as follows how to add twists between each pair of edges so as to obtain the desired embedding of the graph in $\mathbb{R}^3$.

Let $M = \max\{|\lk(C_1, C_2)|\}$
taken over all pairs of disjoint cycles $C_1$ and $C_2$ in $G$.
And let $t = M + \lambda$.

For each $i$, we choose an arc $a_i$ with one endpoint in the interior of $e_i$
and the other endpoint in the interior of $f_i$,
such that $a_i$ is otherwise disjoint from $G$.
In addition, we require that the arcs $a_0, \cdots, a_n$ be pairwise disjoint with neighborhoods that are also pairwise disjoint.  For each $i=0, \cdots, n$,
we change the embedding of $e_i$ and $f_i$ by
adding $3^i t$ positive full twists between them within the neighborhood of $a_i$
(see Figure \ref{twists}).
Note that a given edge will occur in multiple pairs $(e_i,f_i)$ and
hence will be changed within the disjoint neighborhoods of multiple $a_i$'s.

\begin{figure}[h]
\includegraphics[width=3.5in]{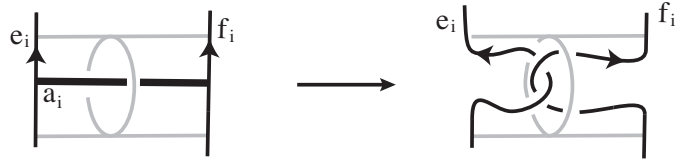}
\caption{We add $3^i t$ positive full twists between $e_i$ and $f_i$ within a
tubular neighborhood of $a_i$}
\label{twists}
\end{figure}

Adding these twists gives us a new embedding $G'$ of $G$.  Fix a pair of disjoint cycles $C_1$ and $C_2$ in $G$,
and let $C'_1$ and $C'_2$ denote the corresponding pair of cycles in $G'$.  Observe that any
crossings between pairs of $a_j$'s or between an $a_j$ and an $e_i$
or $f_i$ that added new crossings between $C_1'$ and $C_2'$
occur in pairs with opposite sign
since we added an integer number of full twists between $e_i$ and $f_i$.
Hence such crossings
do not contribute to the value of $\lk(C'_1, C'_2)$.
Thus
$\lk(C'_1, C'_2) = \lk(C_1, C_2)
+ \epsilon_0 t + \cdots + \epsilon_n 3^n t$,
where $\epsilon_i$ equals 0 or $\pm 1$ depending on
which edges $e_i$ and $f_i$ are in $C_1$ and $C_2$,
and on whether the twists we have added between $e_i$ and $f_i$ induce negative or positive twists between $C_1$ and $C_2$.

Pick the largest $k$ such that $\epsilon_k \ne 0$.
Clearly $k \ge 1$.
Then the pair $(e_k,f_k)$ is lexicographically the largest weighted pair with one edge in $C_1$ and the other in $C_2$.
Therefore the orientations of $C_1$ and $C_2$ are induced by those of $e_k$ and $f_k$.
Since we have added $3^kt$ positive full twists between $e_k$ and $f_k$, it follows that $\epsilon_k = 1$.
Thus we have:
$\lk(C'_1,C'_2)
\geq -M -  (1 +3+ \cdots + 3^{k-1})t +3^k t
= 3^k t  - (3^k -1)t/(3-1)- M
= (3^k +1) t /2 - M
= (3^k +1) (M+\lambda) /2 - M
\geq \lambda$.  Since the choice of the pair of cycles $C_1$ and $C_2$ was arbitrary this proves the proposition.
\end{proof}

\medskip

We saw in Theorem \ref{T:complete} that for any $\lambda>0$ there is a graph $G$ such that every embedding of $G$ contains a link whose linking pattern is $K_n$ with every weight at least $\lambda$.  Recall that the weight of an edge in the linking pattern is defined as the absolute value of the linking number of the associated link in the embedded graph.   It is natural to wonder whether we can remove the absolute value from the conclusion of Theorem \ref{T:complete}.  That is, we would like to know if is there a graph $G$ such that every embedding of $G$ in $\mathbb{R}^3$ contains a link with components $L_1$, \dots, $L_n$ which can be oriented in such a way that $\lk(L_i,L_j)\geq \lambda$ for every $i\not =j$.  We now prove Corollary \ref{C:3link} which shows that the answer to this question is no.

\begin{cor}\label{C:3link}
Let $G$ be a graph.  There is some embedding $G'$ of $G$ in $\mathbb{R}^3$ such that for every three disjoint cycles in $G'$, no matter how they are oriented,
 at least one of the three pairs of oriented cycles will have positive linking number; and there is some embedding $G''$ of $G$ in $\mathbb{R}^3$ such that for every three disjoint cycles in $G''$, no matter how they are oriented,
 at least one of the three pairs of oriented cycles will have negative linking number.
\end{cor}

\begin{proof}  We begin with the embedding $G'$ and the orientation of all of the cycles of $G'$ which is given by Proposition \ref{P:nostrong}.  Thus for any pair of oriented cycles $C_1'$ and $C_2'$ in $G'$, $\lk(C_1',C_2')\geq\lambda>0$.  Now, consider any three disjoint cycles in $G'$.
With the given orientations,
their three pairwise linking numbers are all positive.
If we change the orientation of one of the three cycles, then the linking number between the two unchanged cycles remains positive.
If we change the orientation of two of the three cycles, then
the linking number between the two changed cycles remains positive.
Finally, if we change the orientation of all three cycles, then
all three linking numbers remain positive.  Thus no matter how they are oriented at least one pair will have positive linking number.

Now let $G''$ denote the mirror image of $G'$.
It follows that for every three disjoint cycles in $G''$ and every orientation of these three cycles
at least one of the three will have negative linking number. 

\end{proof}

Thus we have shown that there is no graph whose every embedding contains a link of three cycles where for some orientation of the cycles all the linking numbers are positive, and there is no graph whose every embedding contains a link of three cycles where for some orientation of the cycles all the linking numbers are negative.  By contrast, the next proposition shows that there {\it is} a graph whose every embedding contains a link of three cycles (indeed, of $n$ cycles for any given $n \geq 3$) such that for any orientation of the cycles all the linking numbers have the same sign (either positive or negative).  However, the sign (as well as the cycles) will depend on the specific embedding of the graph in $\mathbb{R}^3$.

\begin{prop} \label{P:Ramsey}
Given any $n \in {\mathbb N}$,
there exists an $r$ such that every embedding of $K_r$ in ${\mathbb R^3}$, with an arbitrary orientation assigned to every cycle,
contains an $n$-component link such that the linking numbers of every pair of components in the oriented link are either
all positive or all negative.
\end{prop}

\begin{proof}
By Ramsey Theory \cite{Ra}, there is an $m$ such that for
any 2-coloring of the edges of $K_m$ with red or blue,
there is a subgraph $K_n$ whose edges are all red or all blue.
We know by Theorem  \ref{T:complete} that, given the above $m$,
every sufficiently large $K_r$ embedded in ${\mathbb R^3}$
contains an $m$-component link whose linking pattern is $K_m$.
Given an embedding of this $K_r$ in $\mathbb{R}^3$, we put an arbitrary orientation on each of the $m$ components of the link in this embedding.
Then we color each edge in the associated linking pattern $K_m$
red or blue according to whether that edge corresponds to
a positive or negative linking number respectively.
Now we know there is a subgraph $K_n$ of this colored $K_m$ whose edges are all red or all blue.
This means that for the $n$-component oriented link corresponding to this monocrhromatic $K_n$ linking pattern, the pairwise linking numbers are either
all positive or all negative.
\end{proof}

\section{Intrinsic knotting of the components}

We extend our definition of weighted linking pattern to include 
knotted components as follows.

\begin{definition}
Given an oriented link $L$ with components $L_1$, \dots, $L_n$, the 
\emph{weighted knotting and linking pattern} of $L$ is the weighted 
linking pattern of $L$ together with a  \emph{weight} assigned to 
each vertex $v_i$ representing the value of $|a_2(L_i)|$.
\end{definition}
\medskip

In light of Theorem \ref{T:complete}, it is natural to ask whether 
for any weighted knotting and linking pattern $\Gamma$, there is a 
graph $G$ such that every embedding of $G$ in $\mathbb{R}^3$ contains 
a link whose weighted knotting and linking pattern is at least as 
``complex'' as $\Gamma$.  The goal of this section is to prove that 
this is indeed the case.  In particular, we will prove the following 
restatement of Theorem \ref{main thm}.

\begin{thm} \label{main thm}
  Let $\alpha \in \mathbb{N}$.  For every $n\in \mathbb{N}$, there is 
a graph $G$ such that every embedding of $G$ in $\mathbb{R}^3$ 
contains a link whose linking pattern is $K_n$ with every edge weight 
and vertex weight at least $\alpha$.\end{thm}
  \medskip

  We begin with some preliminary results.  In the following lemma we 
shall use the notation $A\triangledown\epsilon B$ where $\epsilon\in 
\{0,1\}$.  If $\epsilon=1$, we shall mean $A\triangledown B$.  If 
$\epsilon=0$, we shall mean $A\triangledown \emptyset=A$.

\begin{lemma}\label{old:Lem1} Let $\lambda>0$, and
let $A_1$, $\cdots$, $A_n$,
$B_1$, $\cdots$, $B_{6n+6}$,
be disjoint oriented cycles in a complete graph $K_r$ embedded in 
$\mathbb{R}^3$
such that
$\lk(A_h, B_{i}) \ge \lambda$
for all $h$ and $i$.
Then there exist disjoint cycles
$C_1, C_2, C_3, C_4 \in \{B_i\}$
and an oriented cycle $W'$ in $K_r$ with vertices on
$\bigcup_{i} B_i$
such that $W'$ intersects each of
$C_1, C_2, C_3, C_4$ in exactly one arc with orientation opposite 
that of each $C_i$,
and for every $h$,
$|\lk(A_h, W' \triangledown \epsilon_1 C_1\triangledown\epsilon_2 
C_2\triangledown\epsilon_3 C_3\triangledown\epsilon_4 C_4)| \ge 
\lambda$
for every choice of
$\epsilon_1, \cdots, \epsilon_4 \in \{0,1\}$.

\end{lemma}

\begin{proof}  On each oriented cycle $B_i$, pick adjacent vertices 
$x_i$ and $y_i$ so that the orientation induced on the edge 
$\beta_i=\{x_i,y_i\}$ goes from $x_i$ to $y_i$.  Let $W_0$ be the 
cycle which is the union of
the edges $\beta_i$
and the edges
$\{y_i, x_{i+1}\}$, $1 \le i \le 6n+5$,
and
$\{y_{6n+6}, x_{1}\}$.  Orient $W_0$ so that on each $\beta_i$ the 
orientation goes from $y_i$ to $x_i$.
For the remainder of this proof,
we relabel
$B_1$, $\cdots$, $B_{6n+6}$
as
$B_1^1$, $\cdots$, $B_1^6$,
$B_2^1$, $\cdots$, $B_2^6$, $\cdots$, $B_{n+1}^1$, $\cdots$,
$B_{n+1}^6$, respectively. Figure \ref{W0} illustrates the cycles 
$B_1^1$, $\cdots$, $B_1^6$,
$B_2^1$, $\cdots$, $B_2^6$, $\cdots$, $B_{n+1}^1$, $\cdots$,
$B_{n+1}^6$ together with $W_0$.

\begin{figure}[h]
\includegraphics[width=\textwidth]{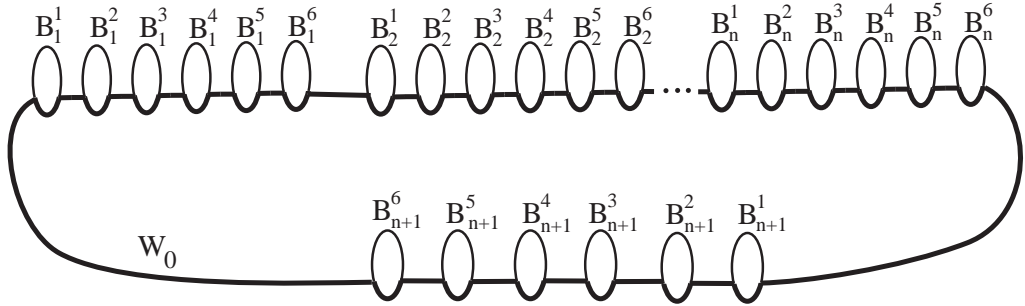}
\caption{$B_1^1$, $\cdots$, $B_1^6$,
$B_2^1$, $\cdots$, $B_2^6$, $\cdots$, $B_{n+1}^1$, $\cdots$,
$B_{n+1}^6$ together with $W_0$.}
\label{W0}
\end{figure}

For $1 \le i \le n$,
let
$W_{i} = W_{i-1} \triangledown B_i^{1} \triangledown  \cdots 
\triangledown B_i^6$.
We have two cases:
\medskip

\emph{Case 1.}
Suppose for some $r \ge 0$,
$\lk(W_r, A_h) \ge 0$ for all $h$.
Then letting
$W' = W_r \triangledown B_{r+1}^6$
and $C_j = B_{r+1}^j$
for $j = 1, 2, 3, 4$
gives us the desired result.
\medskip

\emph{Case 2.}
Suppose for every $k \ge 0$,
$\lk(W_k, A_h) < 0$ for some $h$.
In this case, for each $k$ let $n(k)$ be the number of
$A_h$'s for which
$\lk(W_k, A_h) < 0$.
Then $1 \le n(k) \le n$ for
each $k = 0, 1, \cdots, n$.
Since there are $n+1$ values of $k$, by the Pigeonhole Principle,
$n(r)=n(r')$ for some $r$ and $r'$.  Without loss of generality, $r < r'$.
We see as follows that
$n(k)$ is a non-increasing function of $k$.
For each $h$ and $k$,
$\lk(W_{k}, A_h) > \lk(W_{k-1}, A_h)$
since $\lk(B_k^j, A_h) \ge \lambda > 0$ for all $j$.
It follows that $n(r) = n(r+1)$.
So, for each $h$,
$\lk(A_h,W_r)$ has the same sign as
$\lk(A_h,W_{r+1})$.

Now, as in Case~1, we let
$W' = W_r \triangledown B_{r+1}^6$
and $C_j = B_{r+1}^j$
for $j = 1, 2, 3, 4$.
We verify as follows
that this gives us the desired result.
Fix an $A_h$.
If $\lk(A_h,W_r) \ge 0$,
then clearly
$$\lk(A_h, W' \triangledown \epsilon_1 C_1\triangledown\epsilon_2 
C_2\triangledown\epsilon_3 C_3\triangledown\epsilon_4 C_4) \ge 
\lambda$$
for every choice of
$\epsilon_1, \cdots, \epsilon_4 \in \{0,1\}$,
as desired.
So suppose
$\lk(A_h,W_r) < 0$.
Then
$\lk(A_h,W_{r+1}) < 0$, since $\lk(A_h,W_r)$ has the same sign as
$\lk(A_h,W_{r+1})$.
Now,
for every choice of
$\epsilon_1, \cdots, \epsilon_4 \in \{0,1\}$,
\begin{align*}
\lk(A_h, W'\triangledown \epsilon_1 C_1\triangledown\epsilon_2 
C_2\triangledown\epsilon_3 C_3\triangledown\epsilon_4 C_4) &\le
\lk(A_h, W'\triangledown  C_1\triangledown C_2\triangledown 
C_3\triangledown C_4)\\
&\leq \lk(A_h,W_{r+1}) - \lk(A_h,B_{r+1}^5)\\
&< -\lambda
\end{align*}
as desired.

\end{proof}

\medskip

\begin{lemma}\label{knot}
Let $\lambda>0$, and let
$A_1$, $\cdots$, $A_n$,
$B_1$, $\cdots$, $B_{6n+6}$,
be disjoint oriented cycles in an embedded complete graph $K_r$
such that $\lk(A_h, B_{i}) \ge \lambda$ and $|\lk(B_i, B_{j})| \ge \lambda$
for all $h$, $i$, and $j$.
Then there exists
an oriented cycle $K$ in $K_r$
with vertices on $\bigcup_{i} B_i$
such that
$|a_2(K)| \ge \lambda^2 /16$
and,
for every $h$,
$|\lk(A_h, K)| \ge \lambda$.

\end{lemma}

\begin{proof}
Let $W'$ and $C_1$, $C_2$, $C_3$, and $C_4$ denote the oriented 
cycles given by Lemma \ref{old:Lem1}.  We collapse the four arcs of 
$W'$ which are not in any of the $C_i$ to obtain the pseudograph 
illustrated on the right in Figure \ref{collapseD4}.  We denote this 
pseudograph by $D_4$.

\begin{figure}[h]
\includegraphics{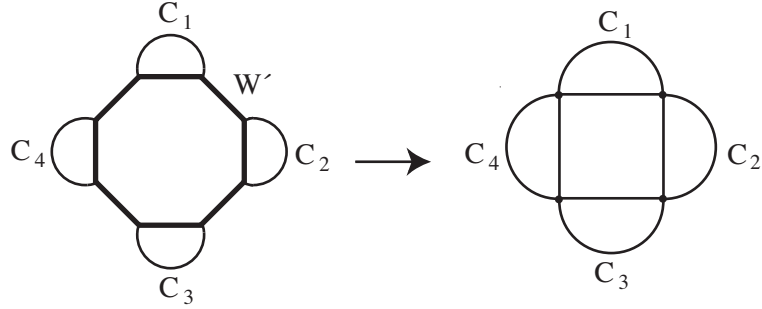}
\caption{We collapse the four arcs of $W'$ not in any $C_i$ to obtain 
the pseudograph $D_4$.}
\label{collapseD4}
\end{figure}
Let $S$ denote the set of all Hamiltonian cycles in $D_4$.  It 
follows from \cite{ty} that

\begin{equation*}  \sum _{Q \in
S  }|a_{2}(Q )|\geq |\lk(C_1,C_3)\lk(C_2,C_4)|
\end{equation*}

Since  $C_1$, $C_2$, $C_3$, $C_4\in \{B_i\}$, we have both $|\lk(C_1, 
C_{3})| \ge \lambda$ and $|\lk(C_2, C_{4})| \ge \lambda$.  Thus $ 
\sum _{Q \in
S  }|a_{2}(Q)|\geq \lambda^2$.  As we are taking the sum of 16 
non-negative integers, there must be some $Q_0\in S $ such
that $ \vert a_{2}(Q_0) \vert \geq
\lambda^2/16$.  Let $K$ denote the cycle of the form $$W' 
\triangledown \epsilon_1 C_1\triangledown\epsilon_2 
C_2\triangledown\epsilon_3 C_3\triangledown\epsilon_4 C_4$$ which 
collapses to $Q_0$ when we collapse the four edges of $W'$ not in any 
$C_i$.  Then $ \vert a_{2}(K) \vert \geq
\lambda^2/16$.  Also, it follows from Lemma \ref{old:Lem1}, that for every $h$,
$|\lk(A_h, K)| \ge \lambda$.  Thus $K$ has both of the required properties.
\end{proof}

\medskip

Roughly speaking,
our goal is to show that
every sufficiently large complete spatial graph
contains a link with a given large number of components,
all with large pairwise linking numbers
and large $a_2$ coefficients.
The idea is that, by Theorem \ref{T:complete},
our large complete spatial graph
contains a link with a large number of components
whose pairwise linking numbers have large absolute value.
We'd like to apply Lemma \ref{knot} repeatedly to this link,
each time increasing the number of components
that have large $a_2$ coefficients,
at the expense of decreasing
the total number of components in the link at hand.
However, to use Lemma \ref{knot},
we need positive linking numbers
between the $A_h$'s and the $B_i$'s.
We accomplish this by giving ourselves the luxury of
starting out with a \emph{lot} of $B_i$'s
and discarding those with negative linking numbers,
as follows.

\begin{lemma}\label{lotsknots} Let $n$, $\lambda \in \mathbb{N}$.
Suppose that a complete graph $K_r$ embedded in $\mathbb{R}^3$
contains an oriented link $L_0$ with $f^n(n)$ components,
where
$f(n) = n - 1+(6n)2^{n-2}$,
such that
the linking number of every pair of components of $L_0$ has
absolute value at least $\lambda$.
Then $K_r$ contains
an oriented link with components $Q_1$, \dots, $Q_n$
such that for every $i \not =j$, $|\lk(Q_i,Q_j)|\geq\lambda$
and $|a_2(Q_i)|\geq\lambda^2/16$.

\end{lemma}

\begin{proof}  Before we begin a recursive argument, we start by 
introducing some variables.  For every $i=1$, \dots, $n$, we let 
$m_{i}=f^{i-1}(n)-1$.  So for each $i$, 
$f^{i}(n)=f(m_{i}+1)=m_{i}+(6m_{i}+6)2^{m_{i}-1}$.  Now for each $i$, 
we let $m_i'=(6m_i+6)2^{m_i-1}$, so that $f^{i}(n)=m_{i}+m_i'$.

We start our recursive argument with the given link $L_0$ which has 
$f^n(n)=m_{n}+m_{n}'$ components.  We
begin by partitioning the components of $L_0$ into two subsets:
$A_1$, $\cdots$, $A_{m_{n}}$,
and
$B_1$, $\cdots$, $B_{m_{n}'}$.
For simplicity,
we shall refer to these two sets
as ``$A$'s'' and ``$B$'s.''  By reversing the orientation of some of 
the $B$'s if necessary,
we can assume they all have positive linking numbers with $A_1$.
Now, $A_2$ has linking numbers of the same sign
with at least half of the $B$'s.
We keep these $B$'s and discard the rest.
By reversing the orientation of $A_2$ if necessary,
we can assume that $A_2$ has positive linking number
with the $B$'s that we kept.
We repeat this process for
$A_3$, $A_4$, \dots, $A_{m_n}$,
each time discarding at most half of the $B$'s.
This reduces the number of $B$'s
by a factor of at most $2^{m_{n}-1}$,
leaving us with at least
${m_{n}'}/{2^{m_{n}-1}}=6m_{n}+6$ remaining
$B$'s, which are each linked to all of the $A$'s with linking number 
at least $\lambda$.

Next we apply Lemma \ref{knot} to the link whose components are 
$A_1$, \dots, $A_{m_n}$ together with $6m_{n}+6$ of the remaining 
$B$'s.  This gives us an oriented knot $Q_1$ which is
linked to all the $A$'s
with absolute value of its linking number
at least $\lambda$
and $|a_2(Q_1)| \ge \lambda^2 /16$.
Let $L_{1}$ be the oriented link whose components are $A_1$, \dots, 
$A_{m_n}$ together with $Q_{1}$.  Then $L_1$ has 
$m_{n}+1=f^{n-1}(n)=m_{n-1}+m_{n-1}'$ components.

We can repeat the above process for the link $L_1$, beginning by 
partitioning the components of $L_1$ into a set of $m_{n-1}$ 
components which includes the knot $Q_1$ and a set of $m_{n-1}'$ 
components.  We abuse notation and refer to the first set as a set of 
$A$'s and the second set as a set of $B$'s.  By applying the above 
argument to these sets of $A$'s and $B$'s we get an oriented link 
$L_{2}$ whose components are the $m_{n-1}$ $A$'s including $Q_1$, 
plus a new oriented knot $Q_{2}$ with $|a_2(Q_2)| \ge \lambda^2 /16$ 
such that $Q_2$ is
linked to all of the $A$'s
with absolute value of its linking number
at least $\lambda$.

We repeat the above process a total of $n-1$ times making sure that 
at each stage all of the $Q_i$'s that we have constructed so far are 
included among the new $A$'s.  In this way we get the desired link 
with components $Q_1$, \dots, $Q_n$, all of whose pairwise linking 
numbers
have absolute value at least $\lambda$
and each component $Q_i$ satisfies
$|a_2(Q_i)| \ge \lambda^2 /16$.

\end{proof}

\medskip

We now prove our main result using
Lemma \ref{lotsknots} and Theorem \ref{T:complete}.

\setcounter{thm}{1}
\begin{thm} \label{main thm} For every $n$, $\alpha\in \mathbb{N}$, 
there is a complete graph $K_r$ such that every embedding of $K_r$ in 
$\mathbb{R}^3$ contains
an oriented link with components $Q_1$, \dots, $Q_n$
such that for every $i\not =j$, $|\lk(Q_i,Q_j)|\geq\alpha$ and 
$|a_2(Q_i)|\geq\alpha$.
\end{thm}

\begin{proof}  Let $\lambda=\mathrm{Max}\{\alpha, 4\sqrt{\alpha}\}$. 
Let $f(n) = n - 1+(6n)2^{n-2}$ and $m=f^n(n)$.  It follows from 
Theorem \ref{T:complete}, that there is a complete graph $K_r$ such 
that every embedding of $K_r$ in $\mathbb{R}^3$
contains an oriented link $L_0$ with $m$ components such that the 
linking number of every pair of components of $L_0$ has absolute 
value at least $\lambda$.  Now it follows from Lemma \ref{lotsknots} 
that every embedding of $K_r$ in $\mathbb{R}^3$ contains
an oriented link with components $Q_1$, \dots, $Q_n$
such that for every $i\not =j$, $|\lk(Q_i,Q_j)|\geq\lambda\geq \alpha$
and $|a_2(Q_i)|\geq\lambda^2/16\geq \alpha$.  \end{proof}

\bigskip

Recall that it follows from Corollary \ref{C:3link} that there is no graph $G$ with the property that every embedding of $G$ in $\mathbb{R}^3$ contains a link with components $Q_1$, \dots, $Q_n$ such that for some orientation of the components, $\lk(Q_i,Q_j)\geq \lambda$ for every $i\not =j$.  Thus the conclusion of Theorem \ref{main thm} cannot be strengthened by removing the absolute value on the linking number.

\small

\normalsize

\end{document}